\documentclass[12pt]{article}

\usepackage[colorlinks=true,citecolor=black,linkcolor=black,urlcolor=blue]{hyperref}
\usepackage{amsthm,amsmath,amssymb}
\usepackage{graphics}
\usepackage{array,caption}

\theoremstyle{plain}
\newtheorem{theorem}{Theorem}
\newtheorem{lemma}[theorem]{Lemma}

\theoremstyle{definition}
\newtheorem{definition}[theorem]{Definition}
\newtheorem{example}[theorem]{Example}

\newtheorem{cons}[theorem]{Construction}

\theoremstyle{remark}

\title{\bf An Easy-to-implement Construction for $(k,n)$-threshold Progressive Visual Secret Sharing Schemes}

\author{Hong-Bin Chen\footnote{Department of Applied Mathematics, National Chung Hsing University, Taichung 40249, Taiwan. Supported by the grant MOST 107-2115-M-035-003-MY2. {\tt Email: hbchen@dragon.nchu.edu.tw}}, Hsiang-Chun Hsu\footnote{Department of Mathematics, Tamkang University, New Taipei City 25137, Taiwan. Supported by the grant MOST 108-2115-M-032-007-MY2. {\tt Email: hchsu0222@gmail.com}} and Justie Su-Tzu Juan\footnote{Corresponding author. Department of Computer Science and Information Engineering,  Nantou 54561, Taiwan. Supported by the grant MOST 108-2221-E-260-008- {\tt Email: jsjuan@ncnu.edu.tw}}
}

\begin{document}

\maketitle


\begin{abstract}
Visual cryptography encrypts the secret image into $n$ shares (transparency) so that only stacking a qualified number of shares  can recover the secret image by the human visual system while no information can be revealed without a large enough number of shares. This paper investigates the $(k,n)$-threshold Visual Secret Sharing (VSS) model, where one can decrypt the original image by stacking at least $k$ shares and get nothing with less than $k$ shares. There are two main approaches in the literature: codebook-based schemes and random-grid-based schemes; the former is the case of this paper. In general, given any positive integers $k$ and $n$, it is not easy to design a valid scheme for the  $(k,n)$-threshold VSS model. In this paper, we propose a simple strategy to construct an efficient scheme for the $(k,n)$-threshold VSS model  for any positive integers $2\leq k\leq n$. The crucial idea is to establish a seemingly unrelated connection between the $(k,n)$-threshold VSS scheme and a mathematical structure -- the generalized Pascal's triangle. This paper improves and extends previous results in four aspects:
\begin{itemize}
\item Our construction offers a unified viewpoint and covers several known results;
\item The resulting scheme has a progressive-viewing property that means the more shares being stacked together the clearer the secret image would be revealed.
\item The proposed scheme can be constructed explicitly and efficiently based on the generalized Pascal's triangle without a computer.
\item Performance of the proposed scheme is comparable with known results.
\end{itemize}

\noindent \textbf{Keywords:} Visual secret sharing scheme, Generalized Pascal's triangle

\end{abstract}

\section{Introduction}

The concept of $(k,n)$-threshold Visual Secret Sharing (VSS) scheme proposed by Naor and Shamir \cite{NS} is to protect a secret image by encoding it to $n$ image shares for a group of $n$ participants so that only $k$ or more participants can recover the secret image. In a VSS scheme, every  participant receives a share consisting of xeroxed transparencies which are stacked to recover the image by the human visual system.
A feature in visual cryptography is that everything is done by the human visual ability and no computation is required in decoding. Therefore, the contrast of the reconstructed images is an important measurement of performance for VSS schemes.

Throughout this paper, the image considered is black-and-white. In the classic model, a VSS scheme expands a pixel $p$ of a binary secret image $P$ into $n$ shares of $m$ sub-pixels each to disguise the pixel $p$ in a way that $p$ can only be visually perceivable with a loss of contrast by stacking at least $k$ transparencies.
According to \cite{NS} by Naor and Shamir, a solution to the $(k,n)$-threshold VSS scheme consists
of two collections of $n \times m$ Boolean matrices $\mathcal{C}_0$ and $\mathcal{C}_1$, called a {\it codebook}. Usually, each of the two collections consists of all matrices of permutating columns of a single matrix, which will be called a {\it basis matrix}. To share a white (black) pixel,
the dealer randomly chooses one of the matrices in $\mathcal{C}_0$ ($\mathcal{C}_1$) according to the codebook, respectively. The chosen matrix
defines the color of the $m$ sub-pixels in each one of the $n$ transparencies. Following Naor and Shamir's idea, several VSS schemes based on the codebook method have been proposed \cite{A1996,CJ,Droste,FL,FLL,HKS2000,ito,NS,Yang}.

Two main disadvantages emerge from the above 1-to-$m$ encoding: the pixel expansion and the image distortion.
To overcome the two potential disadvantages, new schemes with 1-to-1 pixel encoding were proposed \cite{ito,Yang}. The essence of their clever idea is to use the frequency of white pixels in the black and white areas of the recovered image, instead of the contrast of every single pixel, for interpreting black and white pixels by human visual system. To this aim, instead of encoding each pixel into pixels in a whole row, they  encode each pixel of the secret image into only one corresponding pixel in the two basis matrices restricted to a single column (chosen randomly).  In fact, according to \cite{ito,Yang}, one can easily translate a traditional VSS scheme with $m$ pixel expansion to a new VSS scheme using 1-to-1 encoding so as to avoid the expansion and the distortion. Furthermore, the codebook in a conventional scheme still offers the basis matrices $C_0$ and $C_1$ in a new scheme.

Since \cite{NS}, ``how to construct the two basis matrices $C_0$ and $C_1$?'' has been the focus of the study of $(k,n)$-threshold VSS schemes. Meanwhile, {\it progressive visual secret sharing} (PVSS) schemes have been proposed in recent years \cite{CHJ18,CJ,CT,FL,FLL,GLW13,HQ11,LYWC,Shyu07,Shyu09,SC11,WLYWC,WS,YLY,YWN}; namely, the contrast of the proposed scheme is increasing progressively with more and more shares being stacked together. For a particular integer $k$, there are various constructions of basis matrices of a $(k,n)$-threshold PVSS scheme. Hou and Quan \cite{HQ11} found that the codebook of the $(2,n)$-threshold VSS scheme by \cite{NS} can be used as the basis matrices for a $(2,n)$-threshold PVSS scheme. In fact, the codebook of the $(3,n)$-threshold VSS scheme in \cite{NS} is a PVSS scheme. Chen and Juan \cite{CJ} proposed a construction of basis matrices for a $(4,n)$-threshold PVSS scheme.

In general, given any positive integers $k$ and $n$, the codebook of a $(k,n)$-threshold VSS scheme is not easy to generate immediately.
Instead, some turned to find schemes with no codebook, which are based on a so-called {\it random grid method}; see \cite{CT,GLW13,Shyu07,Shyu09,WLYWC,WS,YLY,YWN} for references. Chen and Tsao in 2011 \cite{CT}
obtained a general result for the $(k,n)$-threshold PVSS scheme by using a random grid.
Several studies \cite{GLW13,WS,YWN} have subsequently continued and extended their work. Recently, Yan et al. in 2018 \cite{YLY} have found a new random-grid-based PVSS scheme improving previous results \cite{CT,GLW13,WS} upon visual quality in practice, although the contrast of the proposed
scheme cannot be obtained directly by the given parameters. Random-grid-based VSS schemes can tackle the problem of pixel expansion and codebook design, however, it is thought to suffer from a poor visual quality of the recovered secret image since the background becomes darker when more
shares are stacked together.

Some \cite{Droste,ES,HKS2000,KI,SC11} designed algorithms to find  codebooks for general $(k,n)$-threshold VSS schemes under a specific goal, like optimizing contrast or pixel expansion, where the contrast $\alpha$ (will be defined later) and the pixel expansion $m$ (the size of basis matrices) are the major efficiency standard for a scheme.
 Droste \cite{Droste} was the first to analyze the optimal solution for $\alpha$ and $m$ of a $(k,n)$-threshold VSS scheme, where a method of computing a lower bound of $m$, as well as an upper bound of $\alpha$, was provided. Kotoh and Imai \cite{KI} constructed the basis matrices by solving a linear system and resulted in the same pixel expansion as Droste's construction. With an elaborated linear programming, the optimal pixel expansion of a general $(k,n)$-threshold VSS scheme was obtained by Shyu and Chen \cite{SC11}, and by Eisen and Stinson \cite{ES} for $k=2$. Hofmeister et al. \cite{HKS2000} were the first to compute the exact value of the optimal contrast and proved that the optimal contrast is always achievable by a $(k,n)$-threshold VSS scheme that can be constructed by solving a linear program. However, solving a program can be time-consuming, and thus designing codebooks in an efficient fashion is highly desirable.

A recent paper by Chen, Huang and Juan \cite{CHJ18} has revealed a simple method to design the codebooks of general $(k,n)$-threshold PVSS schemes. Precisely, for any given positive integers $k$ and $n$, the codebook of the $(k,n)$-threshold PVSS scheme can be obtained efficiently and directly by a formula. The general formula therein extends the $(4,n)$-threshold PVSS scheme by Chen and Juan \cite{CJ} to the cases of large $k$. The authors \cite{CHJ18} also show that the performance, including contrast, size constraint and pixel expansion, of their method is comparable with previous results \cite{CJ,CT,FLL,GLW13,NS,WLYWC,WS,YLY,YWN}. However, the formula would be unsatisfied in a theoretical point of view because it seems quite complicated and cannot unify the well-known constructions of $(k,n)$-threshold PVSS schemes for $k=2, 3$ and $n$ by Naor and Shamir \cite{NS}. Over the past few years, all the studies of designing codebooks for the $(k,n)$-threshold PVSS schemes have been worked out case by case. With more designs being revealed, it raises naturally an interesting question that ``is there a single formula to unify all these constructions?'' To this aim, Chen et al. \cite{CHJ18} made a first attempt, which helps us to see the possibility of answering such a question.

The main contribution of this paper is to improve not only the performance but also the design method proposed by Chen et al. \cite{CHJ18}. We show that, with no computer or even with no complicated formula, one can easily construct the codebook of a $(k,n)$-threshold PVSS scheme for any general $2\leq k\leq n$. 
The key to this improvement is to discover a seemingly unrelated but beautiful connection between codebooks of $(k,n)$-threshold PVSS schemes and the well-known generalized Pascal's triangle. All constructions can be obtained explicitly from the generalized Pascal's triangle according to a simple rule. Our results cover several previous known results \cite{CHJ18,CJ,HQ11,NS} as special cases; moreover, using a combinatorial identity, a simple proof is provided in a theoretical way in Combinatorics.

The rest of this paper is organized as follows. In Section 2, some prior knowledge and notations are described. Section 3 introduces the connection between the codebook and the generalized Pascal's triangle.
The main results and proofs are provided in Section 4. Finally, in Section 5, we conclude this paper with some comparison and discussion.

\section{Preliminaries}

This section starts with some definitions and previous results. Throughout this paper matrices are equivalent up to column permutation.
We use $M^n_j$ to denote a binary matrix of size $n\times {n\choose j}$, where every column vector is different from each other and has a constant column weight $j$. For matrices $M_1, M_2, \cdots, M_t$ all have the same row dimension, let $M=[M_0, M_1, \cdots, M_t]$ denote the {\it horizontal concatenation} of these  matrices $M_1, M_2, \cdots, M_t$ along rows. For any positive integer $c$ and any matrix $M$, we define $cM=[M, M, \cdots, M]$ to be concatenating horizontally the matrix $M$ $c$ times.
For example, \[ M^3_0=\begin{bmatrix}
    0        \\
    0     \\
    0
\end{bmatrix}, M^3_1=
\begin{bmatrix}
    1       & 0 & 0 \\
    0       & 1 & 0   \\
    0       & 0 & 1
\end{bmatrix}, M^3_2=
\begin{bmatrix}
    1       & 1 & 0 \\
    1       & 0 & 1   \\
    0       & 1 & 1
\end{bmatrix}, M^3_3=
\begin{bmatrix}
    1      \\
    1  \\
    1
\end{bmatrix},\] and
\[ [2M^3_0, M^3_3]=\begin{bmatrix}
    0 & 0 & 1       \\
    0 & 0 & 1    \\
    0 & 0 & 1
\end{bmatrix}.\]

The matrices $M^n_j$'s are called {\it totally symmetric} matrices in \cite{CJ,HKS2000}. It has been shown in \cite{HKS2000} that the optimal contrast value can always be achieved by a scheme using totally symmetric matrices as its basis matrices. Various constructions for $(k,n)$-threshold VSS schemes based on concatenating some of the matrices $M^n_j$'s are proposed in the literature of Visual Cryptography. To illustrate a simple example, we demonstrate a $(2,2)$-threshold VSS scheme using the following collections of matrices: \\
$\mathcal{C}_0 =\{$all the matrices obtained by permutating the columns of  $\begin{bmatrix}
    0 & 1    \\
    0 & 1
\end{bmatrix}\},$ \\
$\mathcal{C}_1 =\{$all the matrices obtained by permutating the columns of  $\begin{bmatrix}
    0 & 1    \\
    1 & 0
\end{bmatrix}\}.$\\
Traditionally, the dealer first randomly choose two matrices $C_0$ and $C_1$ from the above codebook, respectively. Each pixel of the secret image can be encoded into two sub-pixels of the two shares directly: the first row  represents the share 1 and the second row represents the share 2. No matter what type the secret pixel is, each share receives two sub-pixels with a ``0''(white) and a ``1''(black), and thus cannot obtain enough information on the secret pixel. When two shares are stacked together in a way properly aligning the sub
pixels, we can represent the superimposed share pixels by ``OR'' operation of rows, where ``OR'' operation means $0+0=0$, $0+1=1+0=1$, and $1+1=1$ as shown in the following figure.
\begin{figure}[htbp]
\begin{center}
\scalebox{0.6}{\includegraphics{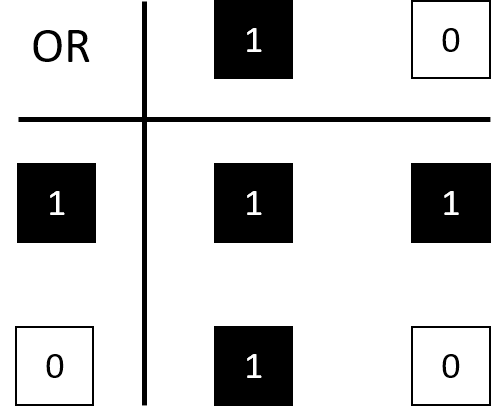}}
\caption{The ``OR'' operation for black and white sub-pixels}
\end{center}
\end{figure}
In the $(2,2)$-threshold VSS scheme above, if the secret pixel is black, then the superimposed share pixels result in two black pixels; however, if the secret pixel is white, then the superimposed share pixels result in one black pixel and one white pixel. This difference makes every secret pixel distinguishable by the human visual system, and thus the secret image is also distinguishable after stacking two shares. Nevertheless, the disadvantage of the pixel expansion could lead to a poor quality of the obtained image.

A method of translating a traditional VSS scheme with $m$ pixel expansion to a new VSS scheme without pixel expansion was proposed in \cite{ito,Yang}.  The following definition is adapted from Yang's definition \cite{Yang} with a slight modification.

\begin{definition}
A $(k,n)$-threshold Visual Secret Sharing (VSS) scheme can be shown as two basis matrices of size $n \times m$: white matrix $C_0$ and black matrix $C_1$.
When sharing a white (resp. black) pixel, the dealer first randomly chooses one column in $C_0$ (resp. $C_1$), and then defines the color level of the pixel in every one of the $n$ shares correspondingly. A $(k,n)$-threshold VSS scheme is considered valid if the following conditions hold:
\begin{itemize}
\item[1.] For any $q$-subset $Q$ of $\{1, 2, \cdots, n\}$ with $q\geq k$, $p_w(Q)-p_b(Q)\geq \alpha >0$, where $p_w(Q)$ and $p_b(Q)$ are the appearance probabilities of 0-entry in the ``OR'' vectors in $C_0$ and $C_1$, respectively, restricted to the $q$ rows indexed by the $q$-subset $Q$, and $\alpha$ is the contrast;
\item[2.] For any $q$-subset $Q$ of $\{1, 2, \cdots, n\}$ with $q<k$, $p_w(Q)=p_b(Q)$.
\end{itemize}
\end{definition}

The first condition defines {\it the contrast $\alpha$}, which is a measurement depending on the appearance probability of the reconstructed white and black pixels; while the second condition is about {\it security} which implies that by stacking less than $k$ shares one cannot obtain enough information on the secret image.

\begin{example}
Take the $(2,2)$-threshold VSS scheme mentioned above for an example. Consider the two basis matrices $C_0 =\begin{bmatrix}
    0 & 1    \\
    0 & 1
\end{bmatrix}$ and
$C_1 =\begin{bmatrix}
    0 & 1    \\
    1 & 0
\end{bmatrix}$. For the only 2-subset $Q=\{1,2\}$, the ``OR'' vector restricted to $Q$, i.e., stacking the two shares, in $C_0$ and $C_1$ are $(0,1)$ and $(1,1)$, respectively. From the definition $p_w(Q)$ and $p_b(Q)$ are the appearance probability of 0-entry in the ``OR'' vector in $C_0$ and $C_1$, respectively, and thus $p_w(Q)-p_b(Q) =1/2 - 0/2=1/2$. Similarly, for any 1-subset $Q=\{1\}$ or $\{2\}$ the ``OR'' vector restrict to $Q$ is indeed a single row and obviously $p_w(Q)-p_b(Q)=1/2 - 1/2=0$. Thus, the two basis matrices are valid for a $(2,2)$-threshold VSS scheme.
\end{example}

Depending on the parameters $k$ and $n$, there are various designs of codebooks for $(k,n)$-threshold VSS schemes based on totally symmetric matrices. Naor and Shamir \cite{NS} proposed the following three constructions:
\begin{itemize}
\item For a solution of a $(2,n)$-threshold PVSS scheme, set the matrices $C_0=[(n-1)M^n_0, M^n_n]$ and $C_1=[M^n_1]$ as the basis matrices.
\item For a solution of a $(3,n)$-threshold PVSS scheme, set the matrices $C_0=[(n-2)M^n_0, M^n_{n-1}]$ and $C_1=[M^n_1, (n-2)M^n_n]$ as the the basis matrices.
\item For a solution of a $(n,n)$-threshold VSS scheme, set the matrices $C_0=[M^n_0, M^n_2, \cdots]$ concatenating all $M^n_j$'s with even weight and $C_1=[M^n_1, M^n_3, \cdots]$ concatenating all $M^n_j$'s with odd weight as the basis matrices.
\end{itemize}
Chen and Juan \cite{CJ} proposed a $(4,n)$-threshold PVSS scheme by setting the matrices $C_0=[\frac{n^2-5n+6}{2}M^n_0, M^n_2, (n-3)M^n_n]$ and $C_1=[(n-3)M^n_1, M^n_{n-1}]$ as the basis matrices.\\

To observe mathematical patterns, we represent the above results by using a single formula of $M^n_j$'s in the form $a_0M^n_0-a_1M^n_1+\cdots +(-1)^ja_jM^n_j+\cdots +(-1)^na_nM^n_n$, where $|(-1)^ja_j|M^n_j$ is in $C_0$ if $(-1)^ja_j\geq 0$, and  $|(-1)^ja_j|M^n_j$ is in $C_1$ for otherwise. Note that the representation is not self-conflicted because the mentioned results are all about the matrices $M^n_j$'s, and moreover each $M^n_j$ appears in at most one of $C_0$ and $C_1$.  We call $(a_0, a_1, \cdots, a_n)$ the {\it coefficient sequence} of the basis matrices.  Previous results are summarized in Table \ref{nscj}. Given a coefficient sequence, we can learn all contents in the basis matrices $C_0$ and $C_1$. Obviously, there is a one-to-one correspondence between the basis matrices and its coefficient sequence. The one-to-one correspondence is crucial in the discussion in the rest of this paper. For the purpose of developing a general construction of basis matrices for $(k,n)$-threshold PVSS schemes, this correspondence suggests that we need not to construct basis matrices directly but only need to simply study its coefficient sequence.

\begin{table}
\begin{center}
    \begin{tabular}{| l | l | l | }
    \hline
    Reference & Model & Coefficient Sequence  \\ \hline
     Naor \& Shamir  & $(2,n)$-threshold & $(n-1, 1, 0, 0, \cdots, 0, (-1)^{n})$
 \\ \hline
    Naor \& Shamir  & $(3,n)$-threshold  & $(n-2, 1, 0, 0, \cdots, 0, (-1)^{n-1}, (-1)^{n-1}(n-2))$  \\ \hline
    Naor \& Shamir  & $(n,n)$-threshold  & $(1, 1, \cdots, 1)$  \\
    \hline
    Chen \& Juan  & $(4,n)$-threshold & $(\frac{n^2-5n+6}{2}, n-3, 1, 0, 0, \cdots, 0, (-1)^{n}, (-1)^{n}(n-3))$ \\
    \hline
    \end{tabular}
    \caption{Previous results with coefficient sequences}\label{nscj}
\end{center}
\end{table}

\section{Connection to Generalized Pascal's triangle}

For $0\leq {\mathcal M}\leq {\mathcal N}$, $\displaystyle {{\mathcal N} \choose {\mathcal M}}$ is commonly used to denote the number of combinations of ${\mathcal N}$ items taken ${\mathcal M}$ at a time, called ${\mathcal N}$ choose ${\mathcal M}$. It is also known as the binomial coefficient of the $x^{{\mathcal M}}$ term in the expression of the binomial power $(1+x)^{{\mathcal N}}$. The well-known Pascal's triangle is a triangular array of the binomial coefficients, where the value of the $({\mathcal N},{\mathcal M})$ entry in the ${\mathcal N}$-th row and ${\mathcal M}$-th column is $\displaystyle {{\mathcal N} \choose {\mathcal M}}$, as shown in the following table.

\begin{center}
\begin{tabular}{>{}l<{\hspace{3pt}}*{13}{c}}
 ${\mathcal N}\setminus {\mathcal M}$  &0&1&2&3&4&5&6&7&8&9&\\
~0 &1&&&&&&&&&&\\
~1 &1&1&&&&&&&&&\\
~2 &1&2&1&&&&&&&&\\
~3 &1&3&3&1&&&&&&&\\
~4 &1&4&6&4&1&&&&&&\\
~5 &1&5&10&10&5&1&&&&&\\
~6 &1&6&15&20&15&6&1&&&&\\
~7 &1&7&21&35&35&21&7&1&&&\\
~8 &1&8&28&56&70&56&28&8&1&&\\
~9 &1&9&36&84&126&126&84&36&9&1&
\end{tabular}
\end{center}
The entries in Pascal's triangle obey a rule that each number is the sum of the two numbers directly above it. This can be written as the well-known Pascal's formula: \[ {{\mathcal N} \choose {\mathcal M}}={{\mathcal N}-1 \choose {\mathcal M}-1}+{{\mathcal N}-1 \choose {\mathcal M}}. \]

In mathematics, Pascal's triangle can be extended to the case of negative row indexes, i.e., ${\mathcal N}<0$. To this aim, define ${{\mathcal N} \choose {\mathcal M}}=1$ for any integer ${\mathcal N}$ if ${\mathcal M}=0$, and then rewrite the above formula as \[ {{\mathcal N}-1 \choose {\mathcal M}}={{\mathcal N} \choose {\mathcal M}}-{{\mathcal N}-1 \choose {\mathcal M}-1}, \] which leads to a simple calculation of the entries for negative row indexes. This extension preserves the property that the values in the ${\mathcal M}$-th column can be viewed as a polynomial function of ${\mathcal N}$ of degree ${\mathcal M}$; namely, \[\displaystyle {{\mathcal N}\choose {\mathcal M}}= \frac{1}{{\mathcal M}!}\prod_{i=1}^{{\mathcal M}}({\mathcal N}+1-i). \]
Such an extension also preserves the property that the values in the ${\mathcal N}$-th row correspond to the coefficients of $(1+x)^{\mathcal N}$ for $|x|<1$. For example, $(1+x)^{-2}= 1-2x+3x^2-4x^3+\cdots$. Obviously, with Pascal's formula, the generalized Pascal's triangle can be produced easily without a computer. Table \ref{pascal} is part of the generalized Pascal's triangle.

\begin{table}[h!]
\begin{center}
\begin{tabular}{>{}l<{\hspace{3pt}}*{13}{c}}
 ${\mathcal N}\setminus {\mathcal M}$  &0&1&2&3&4&5&6&7&8&9&10\\
-8 &1&-8&36&-120&330&-792&1716&-3432&6435&-11440&19448\\
-7 &1&-7&28&-84 &210&-462&924&-1716&3003&-5005&8008\\
-6 &1&-6&21&-56 &126&-252&462&-792&1287&-2002&3003\\
-5 &1&-5&15&-35 &70 &-126&210&-330&495&-715&1001\\
-4 &1&-4&10&-20 &35 &-56 &84 &-120&165&-220&286\\
-3 &1&-3&6 &-10 &15 &-21 &28 &-36 &45 &-55 &66\\
-2 &1&-2&3 &-4  &5  &-6  &7  &-8  &9  &-10 &11\\
-1 &1&-1&1 &-1  &1  &-1  &1  &-1  &1  &-1  &1 \\
~0 &1&0&0&0&0&0&0&0&0&0&0\\
~1 &1&1&0&0&0&0&0&0&0&0&0\\
~2 &1&2&1&0&0&0&0&0&0&0&0\\
~3 &1&3&3&1&0&0&0&0&0&0&0\\
~4 &1&4&6&4&1&0&0&0&0&0&0\\
~5 &1&5&10&10&5&1&0&0&0&0&0\\
~6 &1&6&15&20&15&6&1&0&0&0&0\\
~7 &1&7&21&35&35&21&7&1&0&0&0\\
~8 &1&8&28&56&70&56&28&8&1&0&0\\
~9 &1&9&36&84&126&126&84&36&9&1&0
\end{tabular}
\caption{The generalized Pascal's triangle with small parameters}\label{pascal}
\end{center}
\end{table}

Next, we offer a novel perspective that unifies the results mentioned in Table \ref{nscj}.

\begin{example}\label{ex1}
In the $(2,n)$-threshold VSS scheme by Naor and Shamir \cite{NS}, its coefficient sequence can be found in the generalized Pascal's triangle. Precisely, the sequence starts from the $(n-1,n-2)$ entry and up to the $(-1,n-2)$ entry. Taking $n=4$ for instance, the coefficient sequence $(3,1,0,0,1)$ corresponds to the basis matrices\begin{center} $C_0=[3M^4_0, M^4_4]=\begin{bmatrix}
    0 & 0 & 0&1    \\
    0 & 0 &0&1    \\
    0 & 0 &0&1    \\
    0 & 0 &0&1
\end{bmatrix}$ and $C_1=[M^4_1]= \begin{bmatrix}
    1 & 0 & 0&0    \\
    0 & 1 &0&0    \\
    0 & 0 &1&0    \\
    0 & 0 &0&1
\end{bmatrix}$.
\end{center}
\end{example}

\begin{example}\label{ex2}
According to the $(3,n)$-threshold VSS scheme by Naor and Shamir \cite{NS}, its coefficient sequence starts from the $(n-2,n-3)$ entry and up to the $(-2,n-3)$ entry. Taking $n=4$ for instance, the coefficient sequence $(2,1,0,-1,-2)$ corresponds to the basis matrices \begin{center} $C_0=[2M^4_0, M^4_3]=\begin{bmatrix}
    0 & 0 &1&1&1&0    \\
    0 & 0 &1&1&0&1    \\
    0 & 0 &1&0&1&1    \\
    0 & 0 &0&1&1&1
\end{bmatrix}$ and $C_1=[M^4_1, 2M^4_4]= \begin{bmatrix}
    1 & 0 &0&0&1&1    \\
    0 & 1 &0&0&1&1    \\
    0 & 0 &1&0&1&1    \\
    0 & 0 &0&1&1&1
\end{bmatrix}$.\end{center}
\end{example}

\begin{example}\label{ex3}
According to the $(4,n)$-threshold VSS scheme by Chen and Juan \cite{CJ}, its coefficient sequence starts from the $(n-2,n-4)$ entry and up to the $(-2,n-4)$ entry. Taking $n=5$ for instance, the coefficient sequence $(3,2,1,0,-1,-2)$ corresponds to the basis matrices
 \begin{equation*}
C_0=[3M^5_0, M^5_2, 2M^5_5]=\left[ \begin{array}{@{}*{15}{c}@{}}
  0 &  0 & 0 &1 &1 &1 &1 &0 &0 &0 &0 &0 &0  &1 &1   \\
  0 &  0 & 0 &1 &0 &0 &0 &1 &1 &1 &0 &0 &0  &1 &1   \\
  0 &  0 & 0 &0 &1 &0 &0 &1 &0 &0 &1 &1 &0  &1 &1   \\
  0 &  0 & 0 &0 &0 &1 &0 &0 &1 &0 &1 &0 &1  &1 &1   \\
  0 &  0 & 0 &0 &0 &0 &1 &0 &0 &1 &0 &1 &1  &1 &1
\end{array} \right]
\end{equation*}
and  \begin{equation*}
C_1=[2M^5_1, M^5_4]=\left[ \begin{array}{@{}*{15}{c}@{}}
    1 & 0 &0&0&0 & 1 & 0 &0&0&0 &1 & 1 &1&1&0    \\
    0 & 1 &0&0&0 & 0 & 1 &0&0&0 &1 & 1 &1&0&1   \\
    0 & 0 &1&0&0 & 0 & 0 &1&0&0 &1 & 1 &0&1&1   \\
    0 & 0 &0&1&0 & 0 & 0 &0&1&0 &1 & 0 &1&1&1  \\
    0 & 0 &0&0&1 & 0 & 0 &0&0&1 &0 & 1 &1&1&1
\end{array} \right].
\end{equation*}

\end{example}

\section{Main results}

In this section, we first show the main construction of basis matrices for general $(k,n)$-threshold VSS schemes and then prove its correctness and efficiency.

\begin{cons}\label{con1}
Given any positive integers $k$ and $n$ with $2\leq k\leq n$, let $C_0$ and $C_1$ be the basis matrices whose coefficient sequence starts from the $(n-\lceil\frac{k}{2}\rceil, n-k)$ entry and up to the $(-\lceil\frac{k}{2}\rceil, n-k)$ entry in the generalized Pascal's triangle.
\end{cons}

The proposed design is flexible to apply to any $k\leq n$ and easy to implement because the basis matrices can be obtained immediately from the generalized Pascal's triangle.
It is not difficult to see from Example \ref{ex1}, Example \ref{ex2} and Example \ref{ex3} that each of the mentioned results is a special case of our construction. This construction offers a systematic method to produce basis matrices for $(k,n)$-threshold PVSS schemes.

\begin{theorem}\label{main}
Given any positive integers $k$ and $n$ with $2\leq k\leq n$, the basis matrices $C_0$ and $C_1$ in Construction \ref{con1} are valid for a $(k,n)$-threshold PVSS scheme.
\end{theorem}

To prove Theorem \ref{main}, we need the following extension of combinatorial identities.
Recall that $\displaystyle {{\mathcal N}\choose {\mathcal M}}= \frac{1}{{\mathcal M}!}\prod_{i=1}^{{\mathcal M}}({\mathcal N}+1-i).$
For any formal power series $P(x)$, denote $[x^q]P(x)$ the coefficient of the term $x^q$ in $P(x)$, e.g., if $P(x)=1+2x+5x^3+x^4$, then $[x^3]P(x)=5$.

\begin{lemma}\label{lemma}
Let $s$ and $t$ be nonnegative integers, and let $r$ be an integer. Then
\begin{equation}\label{eq}
\displaystyle \sum_{i=0}^t (-1)^{t-i}{t\choose t-i}{s+r+i\choose s}= \left\{
\begin{array}{llll}
0 & \mbox{ {\rm if} } t\geq s+1,\\
{s+r\choose s-t} & \mbox{ {\rm if} } t\leq s.
\end{array}
\right.
\end{equation}
\end{lemma}
\proof By definition, we have the generating functions \[\displaystyle P(x)=\frac{1}{(1+x)^{s+1}}=\sum_{i\geq 0}{-(s+1)\choose i}x^i=\sum_{i\geq 0}(-1)^i{s+i\choose i}x^i\] and \[ Q(x)=(1+x)^t=\sum_{i\geq 0}{t\choose i}x^i.\] If $r\geq 0$, then \begin{equation*}
\begin{array}{llll}
& & \displaystyle \sum_{i=0}^t (-1)^{t-i}{t\choose t-i}{s+r+i\choose s}& \\
& \displaystyle =&\displaystyle (-1)^{t+r}\sum_{i=0}^t {t\choose t-i}(-1)^{r+i}{s+r+i\choose r+i}& \Big((-1)^{t-i}=(-1)^{t+i}=(-1)^{t+r}\cdot(-1)^{r+i}\Big) \\
& \displaystyle =&\displaystyle (-1)^{t+r}\sum_{i=0}^t \Big( [x^{t-i}]Q(x) \Big)\Big( [x^{r+i}]P(x) \Big)&\\
 &\displaystyle =&\displaystyle (-1)^{t+r}[x^{t+r}]\left((1+x)^t\cdot\frac{1}{(1+x)^{s+1}}\right)& (\mbox{{\rm convolution of the generating functions}})\\
 &\displaystyle =&\displaystyle (-1)^{t+r}[x^{t+r}](1+x)^{t-s-1}& \\
 &\displaystyle =&\displaystyle (-1)^{t+r}{t-s-1\choose t+r}& \\
 &\displaystyle =&\displaystyle (-1)^{t+r}\frac{(t-s-1)(t-s-2)\cdots (t-s-(t+r))}{(t+r)!}.&
\end{array}
\end{equation*}
Therefore, if $r \geq 0$, then we obtain \begin{equation*}
\displaystyle \sum_{i=0}^t (-1)^{t-i}{t\choose t-i}{s+r+i\choose s}= \left\{
\begin{array}{ll}
0 & \mbox{ {\rm if} } t\geq s+1, \\
\displaystyle {s+r\choose s-t} & \mbox{ {\rm if} } t\leq s.
\end{array}
\right.
\end{equation*}
Treat two sides of the above equation as polynomial functions of $r$, i.e., $f(r)$ and $g(r)$. Notice that $f(r)$ and $g(r)$ have a finite degree. As $f(r)=g(r)$ holds for infinitely many $r$ (all integers $r\geq 0$), by the Identity Theorem for polynomials, we have $f(r)=g(r)$ for all integers $r<0$. Hence, Equation (\ref{eq}) holds for any integer $r$, as desired.\qed

\noindent{\bf Proof of Theorem \ref{main}.\quad} Suppose that $C_0$ and $C_1$ are basis matrices corresponding to the coefficient sequence starting from the $(n-\lceil\frac{k}{2}\rceil, n-k)$ entry and up to the $(-\lceil\frac{k}{2}\rceil, n-k)$ entry in the generalized Pascal's triangle. We first verify that $C_0$ and $C_1$ are of the same size $n\times m$; the difference of the numbers of columns in $C_0$ and $C_1$ is
\begin{equation*}\displaystyle \begin{array}{lll}
& \displaystyle \sum_{i=0}^n(-1)^{n-i}{n\choose n-i}{-\lceil\frac{k}{2}\rceil+i\choose n-k} & \\
=&\displaystyle \sum_{i=0}^n(-1)^{n-i}{n\choose n-i}{n-k+r+i\choose n-k} & (\mbox{ {\rm setting } }r=-\lceil\frac{k}{2}\rceil -n+k) \\
=& 0. & (\mbox{ {\rm by Lemma \ref{lemma} with} }t=n>n-k=s)
\end{array}
\end{equation*}

For the conditions concerning security and contrast, we want to prove that if $q<k$ then $p_w(Q)-p_b(Q)=0$ and if $q\geq k$ then $p_w(Q)-p_b(Q)>0$ for any $q$ shares being stacked, where $Q$ is any $q$-subset of $\{1, 2, \cdots, n\}$. Since $C_0$ and $C_1$ are of the same size, it suffices to focus on the difference of the numbers of appearance of 0-entry in the ``or'' vectors in $C_0$ and $C_1$ restricted to a $q$-subset $Q$; that is
\begin{eqnarray}\label{eq2}
\displaystyle \sum_{i=0}^{n-q}(-1)^{n-q-i}{n-q\choose n-q-i}{q-\lceil\frac{k}{2}\rceil+i\choose n-k}
=\sum_{i=0}^{n-q}(-1)^{n-q-i}{n-q\choose n-q-i}{n-k+r+i\choose n-k}
\end{eqnarray}
where the equality holds by setting $r=q-\lceil\frac{k}{2}\rceil -n+k$.

 If $q<k$, then $n-q\geq n-k+1$ and thus by Lemma \ref{lemma} we obtain the above formula (\ref{eq2}) $= 0$. This verifies the security condition. If $q\geq k$, then by Lemma \ref{lemma} we get (\ref{eq2}) $\displaystyle ={n-k+r\choose q-k}
= {q-\lceil\frac{k}{2}\rceil \choose q-k}={q-\lceil\frac{k}{2}\rceil \choose k-\lceil\frac{k}{2}\rceil}$, which is an increasing function of $q$ and has the minimum value 1 when $q=k$. Hence, the contrast $\displaystyle\alpha ={q-\lceil\frac{k}{2}\rceil \choose q-k}/{m} $ is positive and increases progressively for any $q\geq k$ shares being stacked. This completes the proof of Theorem \ref{main}.\qed

We remark that the coefficient sequence in Construction 1 can start from any entry in the $(n-k)$-th column and still works as a valid scheme  because this affects only the parameter $r$ in Lemma \ref{lemma}. Thus, the above proof works as well. Note that the codebook of the $(k,n)$-threshold PVSS scheme proposed in \cite{CHJ18} is a special case whose coefficient sequence starting from the $(n-k+2,n-k)$ entry and up to the $(-k+2,n-k)$ entry. By performing basic calculations, we find that shifting the coefficient sequence of length $n+1$ to the nearly symmetry position in the generalized Pascal's triangle yields a minimum size $m$ of basis matrices. Therefore, starting from the $(n-\lceil\frac{k}{2}\rceil, n-k)$ entry and up to the $(-\lceil\frac{k}{2}\rceil, n-k)$ entry, as in our construction, is the best choice for the  contrast $\alpha$. The coefficient sequences of the mentioned results can be found in the generalized Pascal's triangle, as shown in Figure 2.

\begin{figure}
\begin{center}
\scalebox{0.6}{\includegraphics{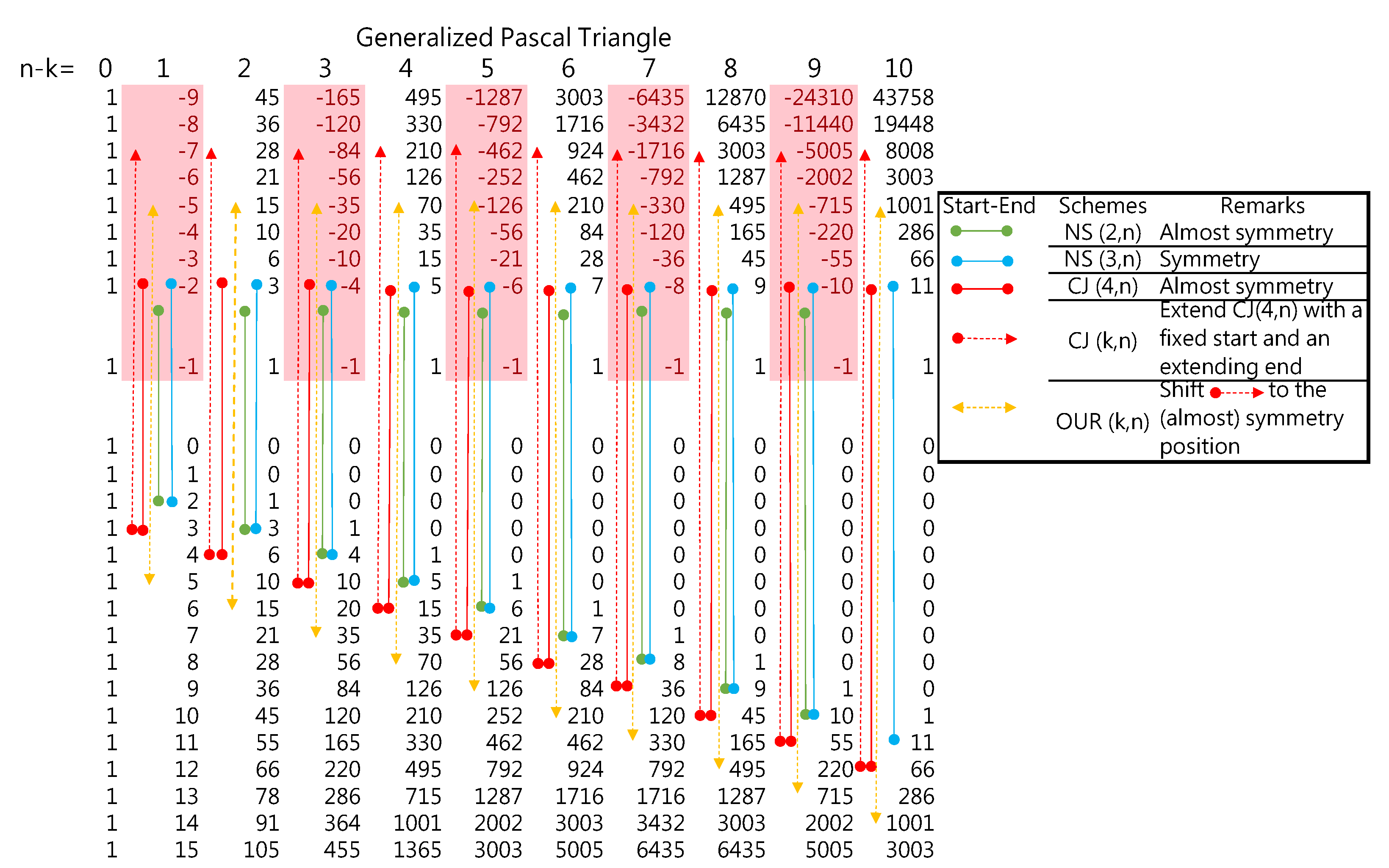}}
\caption{Visualization of the coefficient sequences of the mentioned results in the generalized Pascal's triangle}
\end{center}
\end{figure}

\section{Concluding remarks}
For general $2\leq k\leq n$, we have shown that there exists an explicit construction for the basis matrices $C_0$ and $C_1$ of a $(k,n)$-threshold PVSS scheme. The contrast $\alpha$ of the schemes constructed according to Construction 1 can be described as this formula ${q-\lceil\frac{k}{2}\rceil \choose q-k}/{m}$, where \begin{eqnarray*}
 \displaystyle m &=&\sum_{\substack{i\equiv n ~(\text{mod~} 2)  \\ {-\lceil\frac{k}{2}\rceil+i\choose n-k}\geq 0 }}{n\choose n-i}{-\lceil\frac{k}{2}\rceil+i\choose n-k}+\sum_{\substack{i\equiv n-1 ~(\text{mod~} 2)\\ {-\lceil\frac{k}{2}\rceil+i\choose n-k}< 0 }}{n\choose n-i}{-\lceil\frac{k}{2}\rceil+i\choose n-k}\\
 &=& \frac{1}{2}\sum_i{n\choose n-i}\left|{-\lceil\frac{k}{2}\rceil+i\choose n-k}\right|
\end{eqnarray*} is the number of columns of $C_1$. Some results of the contrast by stacking exactly $k$ shares are listed in  Table \ref{opt}, where theoretically optimal values computed in \cite{HKS2000} are also listed and marked ``OPT''.
\begin{table}
\begin{center}
    \begin{tabular}{ l| ccccccccccccc}
    $k\setminus n$ & 2 & 3 & 4 & 5  & 6  & 7  & 8  & 9   & 10  \\ \hline
                 2 &1/2&1/3&1/4&1/5 &1/6 &1/7 &1/8 &1/9  &1/10  \\
               OPT &1/2&1/3&1/3&3/10&3/10&2/7 &2/7 &5/18 &5/18  \\ \hline
                 3 &   &1/4&1/6&1/8 &1/10&1/12&1/14&1/16 &1/18  \\
               OPT &   &1/4&1/6&1/8 &1/10&1/10&2/21&5/56 &1/12  \\  \hline
                 4 &   &   &1/8&1/15&1/24&1/35&1/48&1/63 &1/80  \\
               OPT &   &   &1/8&1/15&1/18&3/70&3/80&2/63 &1/35
    \end{tabular}\\
    \caption{The contrast $\alpha$ of the resulting $(k,n)$-threshold VSS scheme and the optimal contrast (OPT) by linear programming in \cite{HKS2000}.}\label{opt}
\end{center}
\end{table}
Unfortunately, as shown in Table \ref{opt}, the contrast of the resulting scheme by stacking exactly $k$ shares meets the optimal bound only when $n$ is small. However, the optimal contrast proposed in \cite{HKS2000} counts only for exactly $k$ shares being stacked and their scheme is not a PVSS scheme; thus, by stacking more and more shares, the contrast of our scheme can be higher than their result.

Take $n=8$ and $k=3$ for an example. Their basis matrices are $C_0=[14M^8_0, M^8_6]$ and $C_1=[M^8_2, 14M^8_8]$ while ours are $C_0=[6M^8_0, M^8_7]$ and $C_1=[M^8_4, 6M^8_8]$. It is straightforward to verify that the contrast of our result is better in the cases of $q=7, 8$ shares being stacked. The results of the contrast $\alpha$ by stacking $q\geq 3$ shares are listed in the following Table \ref{HKS} for comparison.
\begin{table}[h!]
\begin{center}
    \begin{tabular}{ c| ccccccccccccc}
                  $q$ && 3  & 4  & 5   & 6   & 7   & 8   \\ \hline
              Ours &   &1/14&2/14&3/14 &4/14 &5/14 &6/14   \\
                   &   &    &    &     &     &$\vee$ &$\vee$   \\
           \cite{HKS2000} &   &4/42&8/42&11/42&13/42&14/42&14/42
    \end{tabular}\\
    \caption{Comparison of the contrast $\alpha$ in $(3,8)$-threshold VSS schemes when stacking $q$ shares}\label{HKS}
\end{center}
\end{table}

The crucial idea behind our construction is all about a well-known mathematical structure -- the generalized Pascal's triangle. 
This provides us a simple way to generalize some of the previous results in \cite{CHJ18,CJ,NS} on $(k,n)$-threshold PVSS schemes systematically.
We also notice that, with small $k$ and $n$, the algorithmic construction based on fully symmetric matrices proposed by Droste \cite{Droste} yields the same basis matrices with ours for a $(k,n)$-threshold PVSS scheme. Intuitively, we suspect, without a proof, that Droste's algorithm eventually outputs the same codebook as that in our construction even if the parameters $k$ and $n$ are large. It would be interesting to investigate the connection behind the coincidence (see Table \ref{droste}). Remarkably, even if such a coincidence exists, our method is much better than that in \cite{Droste} in computation complexity; thus, it does not mean that the main result of this paper is covered by \cite{Droste}. Conversely, such a connection highlights potential applications of the fundamental structures of pure mathematics to the constructions of the codebooks.

\begin{table}[h!]
\begin{center}
    \begin{tabular}{ l| ccccccccccccc}
    $k\setminus n$ & 2 & 3 & 4 & 5  & 6  & 7  & 8  & 9   & 10  \\ \hline
                 2 & 2 & 3 & 4 & 5  & 6  & 7  & 8  & 9   & 10  \\
               3   &   & 4 & 6 & 8  & 10 &12  & 14 & 16  & 18    \\
                 4 &   &   & 8 & 15 &24  &35  &48  & 63  & 80  \\
               5   &   &   &   &16  &30  &48  &70  & 96  &126  \\
                 6 &   &   &   &    &32  &70  &128 &210  &320  \\
               7   &   &   &   &    &    &64  &140 &256  &420  \\
                 8 &   &   &   &    &    &    &128 &315  &640  \\
               9   &   &   &   &    &    &    &    &256  &630  \\
                10 &   &   &   &    &    &    &    &     &512

    \end{tabular}\\
    \caption{The column size $m$ of the resulting $(k,n)$-threshold PVSS scheme is the same with that produced by Droste's algorithm \cite{Droste}.}\label{droste}
    \end{center}
\end{table}

\end{document}